\documentclass[a4,12pt,onecolumn,oneside,notitlepage,final]{article}

\usepackage{latexsym}
\usepackage{amsmath}
\usepackage{amssymb}
\usepackage{amssymb}
\usepackage{color}

\usepackage{cases}
\usepackage{bm}

\date{}
\title{ $ $ \bf Graded extensions in $\bm {K[\mathbb Q,\sigma]$ }}
\begin{document}
\maketitle {$\ $} \\ {Guangming Xie\textsuperscript{1}}, {Miaomiao Wang\textsuperscript{2}}, {Jie Liang\textsuperscript{3}}
\\
{\small School of Mathematics and Statistics, Guangxi Normal
University, Guilin, Guangxi, 541006,  P. R. China\\
 {\textsuperscript{1}gmxie@mailbox.gxnu.edu.cn}\\
 {\textsuperscript{2}1582583007@qq.com}\\
 {\textsuperscript{3}2569646090@qq.com}}
\\[0.5cm]
  \noindent Let $V$ be a total valuation ring of a division ring $K$, $\mathbb Q$ be the additive group of the rational numbers, $Aut(K)$ be the group of automorphisms of $K$. Let $\sigma:\mathbb Q\longrightarrow Aut(K)$ be a group homomorphism and $K[\mathbb{Q},\sigma]$ be the skew group ring of $\mathbb Q$ over $K$. In this paper, we classify graded extensions of $V$ in $K[\mathbb{Q},\sigma]$ into two types and study the structure of them.
\\[0.5cm]
{\small {\it keywords:} Total valuation ring; graded extension; skew group ring.}
\\[0.5cm]
\\[1.0cm]
{\bf 1. Introduction }
\\

 \noindent Let $K$ be a division ring and let $V$ be a total valuation ring of $K$, that is, for any non-zero $k\in K$, either $k\in V$ or $k^{-1}\in V$. We assume that $V\neq K$ throughout this paper. Let $\sigma:\mathbb Q\longrightarrow Aut(K)$ be a group homomorphism and $K[\mathbb{Q},\sigma]$ be the skew group ring of $\mathbb{Q}$ over $K$. For any non-zero $r\in \mathbb Q, r>0$, $K[X^{r},X^{-r};\sigma(r)]$ has a quotient ring $K(X^{r},\sigma(r))$. So we can easily get $K[\mathbb{Q},\sigma]$ has a quotient ring $K(\mathbb{Q},\sigma)$.

 Let $G$ be a non-zero additive subgroup of $\mathbb Q$, $\sigma:G\longrightarrow Aut(K)$ be a group homomorphism and $K[G,\sigma]$ be the skew group ring of $G$ over $K$. Then $K[G,\sigma]$ has a quotient ring $K(G,\sigma)$ by [8, Theorem 2.2]. A graded subring $A=\bigoplus_{g\in G}{A_gX^g}$ of $K[G,\sigma]$ is called a graded total valuation ring of
 $K[G,\sigma]$ if for any non-zero homogeneous element $aX^g$ of $K[G,\sigma]$, either $aX^g\in A$ or $(aX^g)^{-1}\in A$. A graded total valuation ring A of $K[G,\sigma]$ is said to be a graded extension of $V$ in $K[G,\sigma]$ if $A_0=V$.

 A Gauss extension S of $V$ in $K(G,\sigma)$ was defined in [1] as a total valuation ring of $K(G,\sigma)$ with $S\cap K=V$ that satisfies the following condition :

 \begin{center}
 $\alpha S=a_gx^{g}S$
 \end{center}
 for any $\alpha=\sum a_hX^{h}\in K[G,\sigma]$ with $a_gx^{g}S\supseteq a_hx^{h}S$ for all $h$. Then the following results were obtained:\\

 \noindent{\bf Theorem 1.1.} There is a one-to-one correspondence between the set of all Gauss extensions of $V$ in $K(G,\sigma)$ and the set of all graded extensions of $V$ in $K[G,\sigma]$, which is given by $S\longrightarrow S\cap K[G,\sigma]$, where $S$ is a Gauss extension of $V$ in $K(G,\sigma)$([1, (1.8)]). \\

 \noindent{\bf Theorem 1.2.} Let S be a Gauss extension of $V$ in $K(G,\sigma)$ and let $A=S\cap K[G,\sigma]$. Then

(1) The mapping $\varphi:I\longrightarrow I_g=I\cap K[G,\sigma]$ is a one-to-one correspondence between the set of all (right) ideals of $S$ and the set of all graded (right)  ideals of $A$.

(2) $\varphi$ induces a one-to-one correspondence between the set of all prime ideals of $S$ and the set of all graded prime ideals of $A$([1, (2.1)]).\\

 Theorems show that it suffices, in some sense, to study graded extensions in order to study the Gauss extensions. We note that Gauss extensions in [1] were considered in a more general context. Total valuation rings in Ore extensions or in skew polynomial rings have been studied in [2], [3], [4] and [6].

 The structure of graded extensions of $V$ in $K[X,X^{-1};\sigma]\cong K[\mathbb Z,\sigma]$ was studied in [10] and [11]. But for the simplest divisible additive group $\mathbb Q$, the structure of graded extensions of $V$ in $K[\mathbb{Q},\sigma]$ seems more difficult. In [9], we only study graded extensions of $V$ in $K[\mathbb{Q},\sigma]$, which are corresponding some graded maps from $\mathbb Q$ to $\mathbb Z$. In fact, they are just part extensions of
 $V$ in $K[\mathbb{Q},\sigma]$ of Type (e) defined in this paper.

 The aim of this paper is to study the structure of graded extensions of $V$ in $K[\mathbb{Q},\sigma]$.

 In this paper, we classify graded extensions of $V$ in $K[\mathbb{Q},\sigma]$ into two types, that is, Type (I) and Type (II).

 In section 3, we will give a complete description of graded extensions of Type (I).

 In section 4, we will study graded extensions of Type (II) in detail.

 At last, some examples will be given to display some of the various phenomena.\\[1.0cm]
{\bf 2. Preliminaries }\\

\noindent In this section, we collect some notations, definitions and known results. For any additive subgroup $I$ and $J$ of $K$, we set
\begin{center}
 $(J:I)_l=\{a\in K\mid aI\subseteq J\}$,

 $(J:I)_r=\{a\in K\mid Ia\subseteq J\}$,

 $O_l(I)=(I:I)_l$,

 $O_r(I)=(I:I)_r$,

 $I^{-}=\{c^{-1}\mid c\in I,c\neq 0\}$.
\end{center}

\noindent For any left $V$-ideal $I$ and right $V$-ideal $J$, we define
\begin{center}
$^\ast I=\cap \{Wc\mid I\subseteq Wc,c\in K\}$, where $W=O_l(I)$,
$J^\ast=\cap \{cU\mid J\subseteq cU,c\in K\}$, where $U=O_r(J)$.
\end{center}

\noindent For any fixed $g\in\mathbb Q$, $g>0$, we set

\begin{center}
$M_0=V, M_{ig}=A_{g}A_{g}^{\sigma(g)}\cdots A_{g}^{\sigma((i-1)g)}, M_{-ig}=((V:M_{ig})_r)^{\sigma(-ig)}$
\end{center}

\noindent for all $i\in\mathbb N$.

 Let $J(W)$ be the Jacobson radical of $W$. Let $G$ be a non-zero additive subgroup of $\mathbb Q$, $\sigma:G\longrightarrow Aut(K)$ be a group homomorphism and $K[G,\sigma]$ be the skew group ring of $G$ over $K$. $A=\bigoplus_{g\in G}{A_gX^g}$ is a graded extension of $V$ in $K[G,\sigma]$. Let $I$ be a subset of $K$, we will denote $\sigma(g)(I)$ by $I^{\sigma(g)}$ in the following.

  Similar to [5], we have the following definitions.\\

 \noindent {\bf Definition 2.1.} Let $A=\bigoplus_{g\in G}{A_gX^g}$ be a graded extension of $V$ in $K[G,\sigma]$. If for any $g,h\in G$, $g,h>0$, $A_gA_h^{\sigma(g)}=A_{g+h}$, $A_{-g}(A_{-h})^{\sigma(-g)}=A_{-g-h}$, then $A$ is called a graded extension of Type (I).\\

 \noindent {\bf Definition 2.2.} Let $A=\bigoplus_{g\in G}{A_gX^g}$ be a graded extension of $V$ in $K[G,\sigma]$. If $A_gA_h^{\sigma(g)}\neq A_{g+h}$ or $A_{-g}(A_{-h})^{\sigma(-g)}\neq A_{-g-h}$ for some $g,h>0$,  then $A$ is called a graded extension of Type (II).\\

  \noindent {\bf Definition 2.3.} Let $g>0$, $A=\bigoplus_{i\in \mathbb Z}{A_{ig}X^{ig}}$ be a graded extension of $V$ in $K[X^{g},X^{-g};\sigma(g)]$, $O_l(A_g)=W$. We define the following types of graded extensions of $V$ in $K[X^{g},X^{-g};\sigma(g)]$: \\

   Type (a) $A_g=Va=aV^{\sigma(g)}$ and $A_{-g}=V(a^{-1})^{\sigma(-g)}$;

   Type (b) $A_g=Wa\supset aW^{\sigma(g)}$;

   Type (c) $A_g=Wa\subset aW^{\sigma(g)}$;

   Type (d) $A_g=Wa=aW^{\sigma(g)}$, $A_{-g}=J(W)(a^{-1})^{\sigma(-g)}$ and $J(W)^2=J(W)$;

   Type (e) $A_g=Wa=aW^{\sigma(g)}$, $A_{-g}=J(W)(a^{-1})^{\sigma(-g)}$ and $J(W)=Wb^{-1}$ for some $b\in K$;

    Type (f) $^\ast A_g\supset A_g$;

    Type (g) $^\ast A_g=A_g$ and $^\ast M_{ig}$ is not a principal left $W-$ideal for any $i\in \mathbb N$, where $W=O_l(A_g)$;

    Type (h) $^\ast A_g=A_g$ and $^\ast M_{ig}$ is  a principal left $W-$ideal for some $i\in \mathbb N$,  where $W=O_l(A_g)$.\\

   The structure of these types can be found in [10] and [11]. $A$ is a graded extension of Type (I) if and only if $A$ is a graded extension of Type (a), Type (b), Type (c), Type (d), Type (f) or Type (g), and $A$ is a graded extension of Type (II) if and only if $A$ is a graded extension of Type (e) or Type (h) by [5].

   Let $A=\bigoplus_{r\in \mathbb Q}{A_rX^r}$ be a graded extension of $V$ in $K[\mathbb Q,\sigma]$. For any $r\in\mathbb Q$, $r>0$, let $H_r=\bigoplus_{i\in \mathbb Z}{A_{ir}X^{ir}}$. For any $s\in\mathbb Q$, let $W_s=O_l(A_s)$, $W=\bigcup_{s\in \mathbb Q}{W_s}$.\\

 \noindent {\bf Definition 2.4.} Let $A=\bigoplus_{r\in \mathbb Q}{A_rX^r}$ be a graded extension of $V$ in $K[\mathbb Q,\sigma]$, we define the following types of graded extensions in $K[\mathbb Q,\sigma]$:\\

 Type (a)  For any $r>0$, $H_r$ is a graded extension of Type (a) in $K[X^{r},X^{-r};\\
 \sigma(r)]$;

 Type (b)  For any $r>0$, $H_r$ is a graded extension of Type (b) in $K[X^{r},X^{-r};\\\sigma(r)]$;

 Type (c)  For any $r>0$, $H_r$ is a graded extension of Type (c) in $K[X^{r},X^{-r};\\\sigma(r)]$;

 Type (d)   For any $r>0$, $H_r$ is a graded extension of Type (d) in $K[X^{r},X^{-r};\\\sigma(r)]$;

 Type (f)   For any $r>0$, $H_r$ is a graded extension of Type (f) in $K[X^{r},X^{-r};\\\sigma(r)]$;

 Type (g)   For any $r>0$, $H_r$ is a graded extension of Type (g) in $K[X^{r},X^{-r};\\\sigma(r)]$;

 Type (e)  There exist an $r>0$ such that $H_r$ is  a graded extension of Type (II) with $J(W)=b^{-1}W$ for some $b\in K$;

 Type (h)  There exist an $r>0$ such that $H_r$ is a graded extension of Type (II) with $J(W)^2=J(W)$.\\

   Similar to [10, Lemma 1.1], it is easy to get the following lemma.\\

   \noindent {\bf Lemma 2.5.}  Let $G$ be a non-zero additive subgroup of $\mathbb Q$, $A_0=V$. Then $S=\bigoplus_{g\in G}{A_gX^g}\subseteq K[\mathbb Q,\sigma]$ is a graded extension of $V$ in $K[G,\sigma\mid_{G}]$ if and only if

  (i) For any $g\in G$, $A_g$ is an additive subgroup of $K$.

 (ii) For any $g\in G$, $A_g\cup (A_{-g}^{-})^{\sigma(g)}=K$, where $A_{-g}^{-}=\{c^{-1}\mid c\in A_{-g},c\neq 0\}$.

    (iii) For any $g,h\in G$, $A_gA_h^{\sigma(g)}\subseteq A_{g+h}$.\\

   \noindent {\bf Definition 2.6.} Let $G$ be a non-zero additive subgroup of $\mathbb Q$. A map $f:G\longrightarrow \mathbb Z$ is called a graded map if $f(0)=0$, $f(s)+f(t)\leq f(s+t)$ and $f(s)+f(-s)\geq -1$ for all $s,t\in G$.\\

\noindent {\bf Definition 2.7.} Let $r\in\mathbb Q$, $r>0$. A graded map $f:\mathbb Zr\longrightarrow \mathbb Z$ is called a nice map if $f(r)=0$, $f(-r)=-1$.\\

   Let $d$ be a real number, $[d]$ is the Gauss symbol of $d$ and $f_d$, $f_{d}^{(1)}$, $f_{d}^{(-1)}$ are maps from $\mathbb Q$ to $\mathbb Z$. We define $f_d(r)=[rd]$ for all $r\in\mathbb Q$; $f_d^{(1)}(0)=0$, $rd-1\leq f_d^{(1)}(r)<rd$ and $f_d^{(1)}(-r)=-f_d^{(1)}(r)-1$ for all $r\in\mathbb Q$, $r>0$; $f_d^{(-1)}(0)=0$, $f_d^{(-1)}(r)=[rd]$ and $f_d^{(-1)}(-r)=-f_d^{(-1)}(r)-1$ for all $r\in\mathbb Q$, $r>0$. The following was obtained in [9].\\

   \noindent {\bf Theorem 2.8([9]).} Let $\mathbb R$ be the field of real numbers. Then $\{f_d,f_d^{(1)},f_d^{(-1)}\mid d\in\mathbb R\}$ is the set of all graded maps over $\mathbb Q$.\\

    Let $r\in\mathbb Q$, $\alpha_r\in K$, we set

   \begin{center}
   $\alpha_{ir}=\alpha_r(\alpha_r)^{\sigma(r)}\cdots(\alpha_r)^{\sigma((i-1)r)}$ for all $i\in\mathbb N$ and $\alpha_0=1$,

   $\alpha_{-ir}=(\alpha_r^{-1})^{\sigma(-r)}(\alpha_r^{-1})^{\sigma(-2r)}\cdots(\alpha_r^{-1})^{\sigma(-ir)}$ for all $i\in\mathbb N$.
   \end{center}
   \noindent Then $\alpha_{ir}=(\alpha_{-ir}^{-1})^{\sigma(ir)}$ and $(\alpha_{ir})(\alpha_{jr})^{\sigma(ir)}=\alpha_{ir+jr}$ for all $i,j\in\mathbb Z$.

   Let $r\in\mathbb Q$, $A_r\subseteq K$, we set

    \begin{center}
   $M_0=V$, $M_{ir}=A_r(A_r)^{\sigma(r)}\cdots(A_r)^{\sigma((i-1)r)}$,

   $M_{-ir}=((V:M_{ir})_r)^{\sigma(-ir)}$ for all $i\in\mathbb N$.

   \end{center}

   The following two Theorems can be obtained easily from  [10, Theorem 2.11] and [11, Theorem 1.20].\\

   \noindent {\bf Theorem 2.9.} Let $r\in\mathbb Q$, $r>0$, $W$ be an overring of  $V$ and let $H_r=\bigoplus_{i\in\mathbb Z}A_{ir}X^{ir}$ be a subset of $K[X^{r},X^{-r};\sigma(r)]$ with $A_0=V$, $A_r=W\alpha_r=\alpha_rW^{\sigma(r)}$ and $A_{-r}=J(W)\alpha_{-r}$. Suppose that $J(W)=b^{-1}W=Wb^{-1}$ for some $b^{-1}\in J(W)$. Then $H_r$ is a graded extension of $V$ in $K[X^{r},X^{-r};\sigma(r)]$ if and only if the following properties hold:

    (1) There is a nice map $f$ such that $WA_{ir}=Wb^{f(ir)}\alpha_{ir}$ for all $i\in\mathbb Z$.

    (2)(a) If either $W=V$ or $f(ir)+f(-ir)=-1$ for all $i\in\mathbb Z$ with $i\neq 0$, then $A_{ir}=Wb^{f(ir)}\alpha_{ir}$ for all $i\in\mathbb Z$ with $i\neq 0$.

   \hspace{1.2em}(b) If $W\neq V$ and there is an $l\in\mathbb N(l\geq 2)$ with $f(lr)+f(-lr)=0$(assume that $l$ is the smallest natural number for this property), then $A_{ir}=Wb^{f(ir)}\alpha_{ir}$ for all $i\notin l\mathbb Z$ and $B=\bigoplus_{j\in\mathbb Z}A_{jlr}X^{jlr}$ is a graded extension of $V$ in $K[X^{lr},X^{-lr};\sigma(lr)]$ with
   $Wb^{f(jlr)-1}\alpha_{jlr}\subset A_{jlr}\subset Wb^{f(jlr)}\alpha_{jlr}$ for all $j\in\mathbb Z$.

   Furthermore, assume that the above conditions hold. Set $\beta_{jr}=b^{f(jr)}\alpha_{jr}$, then $W\beta_{jr}=\beta_{jr}W^{\sigma(jr)}$ for all $j\in\mathbb Z$.\\

   \noindent {\bf Theorem 2.10.} Let $r\in\mathbb Q$, $r>0$, $H_r=\bigoplus_{i\in\mathbb Z}A_{ir}X^{ir}$ be a subset of $K[X^{r},X^{-r};\sigma(r)]$ such that $A_0=V$, $A_r$ is a left $V$ and right $\sigma(V)$-ideal with $^\ast A_r=A_r$ which is not a principal left $W$-idea, where $W=O_l(A_r)$. Suppose that $M_{lr}=J(W)c_{lr}$ for some $c_{lr}\in K$ and $l\in\mathbb N$(assume that $l$ is the smallest natural number for this property). Then $H_r$ is a graded extension of $V$ in $K[X^{r},X^{-r};\sigma(r)]$ if and only if

  (1) $A_{ir}=M_{ir}$ for all $i\in\mathbb Z\backslash l\mathbb Z$.

  (2) $B=\bigoplus_{j\in\mathbb Z}A_{jlr}X^{jlr}$ is a graded extension of $V$ in $K[X^{lr},X^{-lr};\sigma(lr)]$ with $J(W)c_{lr}\subseteq A_{lr}\subseteq Wc_{lr}$.

   Furthermore, assume that the above conditions hold. Then $Wc_{jlr}=c_{jlr}W^{\sigma(jlr)}$ for all $j\in\mathbb Z$.
   \\[1.0cm]
{\bf 3. Graded extensions of Type (I) in ${\bm {K[\mathbb Q,\sigma ]}}$}\\

\noindent In this section, we will give a complete classification of graded extensions of Type (I) in $K[\mathbb Q,\sigma]$.

In [5], [10] and [11], we studied graded extensions of $V$ in $K[X,X^{-1};\sigma]$ in detail. We will us the results of them.

Let $A=\bigoplus_{r\in \mathbb Q}{A_rX^r}$ be a graded extension of $V$ in $K[\mathbb Q,\sigma ]$. For any $r\in\mathbb Q$, $r>0$, let $H_r=\bigoplus_{i\in \mathbb Z}{A_{ir}X^{ir}}$. By [10] and [11], we can easily get the following lemmas.\\

\noindent {\bf Lemma 3.1.} Let $A=\bigoplus_{r\in \mathbb Q}{A_rX^r}$ be a graded extension of $V$ in $K[\mathbb Q,\sigma ]$. Let $r\in\mathbb Q$, $r>0$. Then $H_r$ is a graded extension of Type (I) if and only if $H_{nr}$ is a graded extension of Type (I) for any $n\in\mathbb N$.\\


  \noindent {\bf Lemma 3.2.} Let $A=\bigoplus_{r\in \mathbb Q}{A_rX^r}$ be a graded extension of $V$ in $K[\mathbb Q,\sigma ]$, $r\in\mathbb Q$, $r>0$. If $H_r$ is a graded extension of Type (a)(or (b), (c), (d), (f), (g)). Then $H_{nr}$ is a graded extension of Type (a)(or (b), (c), (d), (f), (g)) for any $n\in\mathbb N$.\\

  As we know, for Type (e) and Type (h), we can not get similar results.

  The following is the main result of this section.\\

  \noindent {\bf Theorem 3.3.} Let $A=\bigoplus_{r\in \mathbb Q}{A_rX^r}$ be a graded extension of $V$ in $K[\mathbb Q,\sigma ]$. Then the following are equivalent:

   (1) $A$ is a graded extension of Type (I).

   (2) $A$ is a graded extension of Type (a), Type (b), Type (c), Type (d), Type (f) or Type (g).\\

  \noindent {\bf Proof.} Let $A$ be a graded extension of $V$ in $K[\mathbb Q,\sigma]$ of Type (I). For any $r,s\in \mathbb Q$, $r>0,s>0$, $H_r=\bigoplus_{i\in \mathbb Z}{A_{ir}X^{ir}}$ and $H_s=\bigoplus_{i\in \mathbb Z}{A_{is}X^{is}}$ are graded extensions of Type (I). Let $r=n_1/m_1$, $s=n_2/m_2$, $m_1,m_2,n_1,n_2\in\mathbb N$. Set $m=m_1m_2$. Then $H_{1/m}$ is a graded extension of Type (I). By [5, Theorem 2.3], $H_{1/m}$ is a graded extension of Type (a), Type (b), Type (c), Type (d), Type (f) or Type (g). Furthermore, $H_r$, $H_s$, $H_{1/m}$ are graded extensions of the same type by Lemma 3.2. Hence $A$ is a graded extension of Type (a), Type (b), Type (c), Type (d), Type (f) or Type (g).

  Let $A$ be a graded extension of Type (a), Type (b), Type (c), Type (d), Type (f) or Type (g). Let $r$, $s$ and $m$ be the same as the above. Then $H_r\subseteq H_{1/m}$ and $H_s\subseteq H_{1/m}$ and $H_{1/m}$ is a graded extension of Type (a), Type (b), Type (c), Type (d), Type (f) or Type (g). So $A_rA_s^{\sigma(r)}=A_{r+s}$, $A_{-r}(A_{-s})^{\sigma(-r)}=A_{-r-s}$. Therefore $A$ is a graded extension of Type (I).$\hfill\square$\\[1.0cm]
{\bf 4. Graded extensions of Type (II) in ${\bm {K[ \mathbb Q,\sigma ]}}$}\\

\noindent In this section, we will study graded extensions of Type (II) in $K[\mathbb Q,\sigma ]$.

  Let $A=\bigoplus_{r\in \mathbb Q}{A_rX^r}$ be a graded extension of Type (II) in $K[\mathbb Q,\sigma ]$. Let
$W_r=O_l(A_r)$ for any $r\in\mathbb Q$ and $W=\bigcup_{r\in\mathbb Q}W_r$. If $J(W)$ is a principal $W$-ideal, then $A$ is a graded extension of Type (e). If $J(W)$ is not finitely generated $W$-ideal, then $A$ is a graded extension of Type (h). For any $r\in\mathbb Q$, $r>0$, let $H_r=\bigoplus_{i\in \mathbb Z}{A_{ir}X^{ir}}$. The following lemmas can be obtained by [5], [10] and [11].\\

\noindent {\bf Lemma 4.1.} Let $s\in\mathbb Q$, $s>0$. Assume that $H_s$ is a graded extension of Type (II), then for any $n\in\mathbb N$, $H_{s/n}$ is a graded extension of Type (II).\\

\noindent {\bf Lemma 4.2.}  Let $s\in\mathbb Q$, $s>0$. Assume that $H_s=\bigoplus_{i\in \mathbb Z}{A_{is}X^{is}}$ is a graded extension of Type (II) in $K[X^{s},X^{-s};\sigma(s)]$. Then for any $i\in\mathbb Z$, $W_{is}\subseteq W_s=W_{-s}$.\\

By the proof of Theorem 3.3, we can deduce that there is a $t\in\mathbb Q$, $t>0$ such that $H_t=\bigoplus_{i\in \mathbb Z}{A_{it}X^{it}}$ is a graded extension of Type(II) in $K[X^{t},X^{-t};\sigma(t)]$. Let $t=n/{k_0}$ with $k_0, n\in\mathbb N$. Then  $H_{1/k_0}=\bigoplus_{i\in \mathbb Z}{A_{i/{k_0}}X^{i/{k_0}}}$  is a graded extension of Type (II) by lemma 4.1. So we can get the following lemma.\\

\noindent {\bf Lemma 4.3.} Let $A=\bigoplus_{r\in \mathbb Q}{A_rX^r}$ be a graded extension of Type (II) in $K[\mathbb Q,\sigma ]$. Then there exist a ${k_0}\in \mathbb N$ such that $H_{1/{k_0}}=\bigoplus_{i\in \mathbb Z}{A_{i/{k_0}}X^{i/{k_0}}}$  is a graded extension of Type (II) in $K[X^{1/{k_0}},X^{-1/{k_0}};\sigma(1/{k_0})]$.\\

In the following, we always assume that $H_{1/{k_0}}$ is a graded extension of Type (II).

In this section, we always assume that there is a $t\in\mathbb Q$ with $O_l(A_t)=W$.

At first, we study graded extensions of Type (e), that is, $J(W)=b^{-1}W$ for some $b\in K$. We have the following lemmas.\\

\noindent {\bf Lemma 4.4.} Let $A=\bigoplus_{r\in \mathbb Q}{A_rX^r}$ be a graded extension of Type (e) in $K[\mathbb Q,\sigma ]$. Assume that there exist a $t\in\mathbb Q$ with $O_l(A_t)=W$. Then there exist an $l\in\mathbb N$ such that $H_{1/l}=\bigoplus_{i\in \mathbb Z}{A_{i/l}X^{i/l}}$ is a graded extension of $V$ in $K[X^{1/l},X^{-1/l};\sigma(1/l)]$ of Type (e) with $O_l(A_{1/l})=W$.\\

\noindent {\bf Proof.} We can assume that $t\neq 0$. Let $t=n/m$ with $m\in\mathbb N$, $n\in\mathbb Z$. By Lemma 4.3, there exist a $k_0\in\mathbb N$ such that $H_{1/k_0}$ is a graded extension of Type (II). Let $l={k_0}m$. Then $H_{1/l}$ is a graded extension of Type (II) by Lemma 4.1. $W\supseteq W_{1/l}\supseteq W_t=W$ by Lemma 4.2. Since $J(W)=b^{-1}W$ for some $b\in K$, $H_{1/l}$ is a graded extension of Type (e).$\hfill\square$\\

To study the graded extensions of Type (e) in $K[\mathbb Q,\sigma]$, nice maps are very crucial.
Let $r\in\mathbb Q$, $r>0$, $H_r=\bigoplus_{i\in \mathbb Z}{A_{ir}X^{ir}}$ be a graded extension of $V$ in $K[X^{r},X^{-r};\sigma(r)]$ of Type (e) with $O_l(A_r)=W$. Then there are $\alpha_{ir}\in K$ for every $i\in\mathbb Z$ and a nice map $f$ with $WA_{ir}=Wb^{f(ir)}\alpha_{ir}$. Set $\beta_{ir}=b^{f(ir)}\alpha_{ir}$. Then we can easily get the following lemma.\\

\noindent {\bf Lemma 4.5.} Let $r\in\mathbb Q$, $r>0$, $H_r=\bigoplus_{i\in \mathbb Z}{A_{ir}X^{ir}}$ be a graded extension of $V$ in $K[X^{r},X^{-r};\sigma(r)]$ of Type (e) with $O_l(A_r)=W$. Then we have the following:

 (i) For any $i\in\mathbb Z$, there is a $\beta_{ir}$ with $WA_{ir}=W\beta_{ir}=\beta_{ir}W^{\sigma(ir)}$.

 (ii) $D_r=\bigoplus_{i\in \mathbb Z}{(WA_{ir})X^{ir}}$ is a graded extension of $W$ in $K[X^{r},X^{-r};\sigma(r)]$ of Type (e).\\

The following lemma can be obtained easily by induction.\\

\noindent {\bf Lemma 4.6.} Let $r\in\mathbb Q^{+}$, $f:\mathbb Zr\longrightarrow\mathbb Z$ be a graded map. Then for any $n\in\mathbb N$. $f(nr)<n(f(r)+1)$.\\

\noindent {\bf Lemma 4.7.}  Let $A=\bigoplus_{r\in \mathbb Q}{A_rX^r}$ be a graded extension of Type (e) in $K[\mathbb Q,\sigma ]$. Assume that there exist a $t\in\mathbb Q$ with $O_l(A_t)=W$. Then for any $r\in\mathbb Q$, there is a $\beta_r\in K$ with $WA_r=W\beta_r=\beta_rW^{\sigma(r)}$.\\

\noindent {\bf Proof.} Let $r=n_1/m_1$, $m_1\in\mathbb N$, $n_1\in\mathbb Z$. By Lemma 4.4, there is an $l\in\mathbb N$ such that $H_{1/l}$ is a graded extension of Type (e) with $W_{1/l}=W$. Set $n=lm_1$. Then $H_{1/n}$ is a graded extension of Type (e) with $W_{1/n}=W$. Since $r=ln_1/n$, $WA_r=W\beta_r=\beta_rW^{\sigma(r)}$ by Lemma 4.5.$\hfill\square$\\

Let $A=\bigoplus_{r\in \mathbb Q}{A_rX^r}$ be a graded extension of $V$ in $K[\mathbb Q,\sigma ]$ of Type (e). For any $r\in\mathbb Q$, set $C_r=A_{-r}$, $Y^r=X^{-r}$, $\tau(r)=\sigma(-r)$. By [11, Proposition 1.5], $C=\bigoplus_{r\in \mathbb Q}{C_rY^r}$ is a graded extension of $V$ in $K[\mathbb Q,\tau]$. For any $r\in\mathbb Q$, set $B_r=WA_r$, $E_r=WC_r$, then $E_r=WC_r=WA_{-r}=B_{-r}$. By Lemma 4.7, $B=\bigoplus_{r\in \mathbb Q}{B_rX^r}=\bigoplus_{r\in \mathbb Q}{(WA_r)X^r}$ is a graded extension of $W$ in $K[\mathbb Q,\sigma ]$ and $E=\bigoplus_{r\in\mathbb Q}E_rY^r=\bigoplus_{r\in\mathbb Q}(WC_r)Y^r$ is a graded extension of $W$ in $K[\mathbb Q,\tau]$. Then by Lemma 4.4 and 4.5, there is an $l\in\mathbb N$, both $\bigoplus_{i\in\mathbb Z}B_{i/l}X^{i/l}$ and $\bigoplus_{i\in\mathbb Z}E_{i/l}X^{i/l}$ are graded extensions of Type (e). Assume that $r=n_1/m_1$ for some $m_1\in\mathbb N$, $n_1\in\mathbb Z$, $(n_1,m_1)=1$. Set $q=lm_1$. For any $n\in\mathbb N$, $B_{1/(nq)}=WA_{1/{(nq)}}=W\beta_{1/{(nq)}}$, $E_{1/{(nq)}}=WA_{-1/{(nq)}}=W\beta_{-1/{(nq)}}$. Set $\gamma_{(n,0,q)}=1$, $\delta_{(-n,0,q)}=1$. For any $i\in\mathbb N$,

 \begin{center}
 $\gamma_{(n,i,q)}=\beta_{\frac{1}{nq}}\beta_{\frac{1}{nq}}^{\sigma (\frac{1}{nq})} \cdot\cdot\cdot \beta_{\frac{1}{nq}}^{\sigma((i-1)\cdot\frac{1}{nq})}, \gamma_{(n,-i,q)}=(\gamma_{(n,i,q)}^{-1})^{\sigma(-\frac{i}{nq})}$,

 $\delta_{(-n,i,q)}=\beta_{-\frac{1}{nq}}(\beta_{-\frac{1}{nq}})^{\tau(\frac{1}{nq})} \cdot\cdot\cdot (\beta_{-\frac{1}{nq}})^{\tau((i-1)\cdot\frac{1}{nq})}, \delta_{(-n,-i,q)}=(\delta_{(-n,i,q)}^{-1})^{\tau(-\frac{i}{nq})}$.
\end{center}

Since  $\bigoplus_{i\in \mathbb Z}{B_{\frac{i}{nq}}X^{\frac{i}{nq}}}$ is a graded extension of Type (e), there is a nice map $g_{\frac{1}{nq}}:\frac{1}{nq}\mathbb Z\longrightarrow \mathbb Z$ such that

\begin{center}
 $WA_{\frac{1}{nq}}=W\beta_{\frac{1}{nq}}=W\gamma_{(n,1,q)}$, $WA_{-\frac{1}{nq}}=W\beta_{-\frac{1}{nq}}=Wb^{-1}\gamma_{(n,-1,q)}$.
\end{center}

\noindent For any $i\in\mathbb Z$, $B_{\frac{i}{nq}}=Wb^{g_{\frac{1}{nq}}(\frac{i}{nq})}\gamma_{(n,i,q)}$.

For any $n\in\mathbb N$, $\bigoplus_{i\in \mathbb Z}{E_{\frac{i}{nq}}Y^{\frac{i}{nq}}}$ is a graded extension of Type (e), then there is a nice map $h_{\frac{1}{nq}}:\frac{1}{nq}\mathbb Z\longrightarrow \mathbb Z$ such that

\begin{center}
 $WC_{\frac{1}{nq}}=WA_{-\frac{1}{nq}}=W\beta_{-\frac{1}{nq}}=W\delta_{(-n,1,q)}$, $WC_{-\frac{1}{nq}}=WA_{\frac{1}{nq}}=W\beta_{\frac{1}{nq}}=Wb^{-1}\delta_{(-n,-1,q)}$.
\end{center}

\noindent For any $i\in\mathbb Z$, $E_{\frac{i}{nq}}=Wb^{h_{\frac{1}{nq}}(\frac{i}{nq})}\delta_{(-n,i,q)}$.

Note that $r=\frac{n_1}{m_1}=\frac{n_1l}{m_1l}=\frac{n_1l}{q}$, $g_{\frac{1}{nq}}(r)=g_{\frac{1}{nq}}(\frac{nn_1l}{nq})$, $h_{\frac{1}{nq}}(r)=h_{\frac{1}{nq}}(\frac{nn_1l}{nq})$. Let \begin{center}
$l_r=sup\{g_{\frac{1}{nq}}(r)\mid n\in\mathbb N\}$,

$k_r=sup\{h_{\frac{1}{nq}}(r)\mid n\in\mathbb N\}$.\\
\end{center}
\noindent {\bf Remark 4.8.} By the definition of $l_r$ and $k_r$, we can easily get that $l_r$ and $k_r$ are not depend
on the choice of $l$.\\

Now we consider the boundedness of $l_r$ and $k_r$.\\

\noindent {\bf Lemma 4.9.} For any $r_1,r_2\in\mathbb Q^{+}$, $l_{r_1}<+\infty \Longleftrightarrow l_{r_2}<+\infty$, $k_{r_1}<+\infty \Longleftrightarrow k_{r_2}<+\infty$.\\

\noindent {\bf Proof.} Let $r_1=\frac{n_1}{m_1}$, $r_2=\frac{n_2}{m_2}$, $m_1,m_2,n_1,n_2\in\mathbb N$, $(m_1,n_1)=1$, $(m_2,n_2)=1$.

\begin{center}
 $g_{\frac{1}{nq_1}}(r_1)=g_{\frac{1}{nq_1}}(\frac{nn_1l}{nq_1})\leq g_{\frac{1}{(m_2n)q_1}}(\frac{m_2nn_1l}{m_2nq_1})=g_{\frac{1}{(m_2n)q_1}}(r_1)$,
 $g_{\frac{1}{nq_2}}(r_2)=g_{\frac{1}{nq_2}}(\frac{nn_2l}{nq_2})\leq g_{\frac{1}{(m_1n)q_2}}(\frac{m_1nn_2l}{m_1nq_2})=g_{\frac{1}{(m_1n)q_2}}(r_2)$,
\end{center}

\noindent since $q_1=lm_1$, $q_2=lm_2$, $m_2nq_1=lm_1m_2n=m_1nq_2$. If $l_{r_1}<+\infty$. Assume $l_{r_1}=M$. Then

 \begin{center}
$g_{\frac{1}{m_1nq_2}}(\frac{nl}{m_1nq_2})\leq g_{\frac{1}{m_1nq_2}}(\frac{m_2n_1nl}{m_1nq_2})\leq M$.
\end{center}

\noindent By Lemma 4.6, we obtain

\begin{equation}
\begin{aligned}
   g_{\frac{1}{m_1nq_2}}(r_2)&=g_{\frac{1}{m_1nq_2}}((m_1n_2)\cdot{\frac{nl}{m_1nq_2})}\\
  &< m_1n_2(g_{\frac{1}{m_1nq_2}}(\frac{nl}{m_1nq_2})+1)\\
  &\leq m_1n_2(M+1).
  \nonumber
\end{aligned}
\end{equation}

\noindent So $l_{r_2}<+\infty$. By symmetry, we have $l_{r_1}<+\infty \Longleftrightarrow l_{r_2}<+\infty$. Similarly, we can get $k_{r_1}<+\infty \Longleftrightarrow k_{r_2}<+\infty$.$\hfill\square$\\

\noindent {\bf Lemma 4.10.} Let $r\in\mathbb Q^{+}$ and assume that $l_r<+\infty$, then $k_r=+\infty$.\\

\noindent {\bf Proof.} Let $r=\frac{n_1}{m_1}$, $m_1,n_1\in\mathbb N$, $(m_1,n_1)=1$, $q=lm_1$. Then $\bigoplus_{i\in\mathbb Z}B_{\frac{i}{nq}}X_{\frac{i}{nq}}$ is a graded extension of Type (e) in
$K[X^{\frac{1}{nq}},X^{-\frac{1}{nq}};\sigma(\frac{1}{nq})]$. Assume $l_r=m\in\mathbb N$. 
For any big enough $M\in\mathbb N$, let $n=M+m+2$. Since $\bigoplus_{i\in\mathbb Z}E_{\frac{i}{nq}}Y_{\frac{i}{nq}}$ is a graded extension of Type (e),

 \begin{center}
$B_{-\frac{1}{nq}}=W\beta_{-\frac{1}{nq}}=W\delta_{(-n,1,q)}=Wb^{-1}\gamma_{(n,-1,q)}$,
\end{center}
\begin{equation}
\begin{aligned}
   W\delta_{(-n,nn_1l,q)}&=W\beta_{-\frac{1}{nq}}(\beta_{-\frac{1}{nq}})^{\sigma(-\frac{1}{nq})}\cdots(\beta_{-\frac{1}{nq}})^{\sigma((nn_1l-1)\cdot(-\frac{1}{nq}))}\\
  &= Wb^{-nn_1l}\gamma_{(n,-nn_1l,q)}.
   \nonumber
\end{aligned}
\end{equation}

\noindent Then

\begin{center}
$E_r=Wb^{h_{\frac{1}{nq}}(r)}\delta_{(-n,nn_1l,q)}=Wb^{{h_{\frac{1}{nq}}(r)}-nn_1l}\gamma_{(n,-nn_1l,q)}$.
\end{center}

\noindent Furthermore, $E_r=B_{-r}=Wb^{g_{\frac{1}{nq}}(-r)}\gamma{(n,-nn_1l,q)}$. Hence
\begin{center}
${h_{\frac{1}{nq}}(r)}-nn_1l=g_{\frac{1}{nq}}(-r)\geq -g_{\frac{1}{nq}}(r)-1\geq -m-1$.
\end{center}
\begin{equation}
\begin{aligned}
   {h_{\frac{1}{nq}}(r)}&\geq nn_1l-m-1\\
  &\geq n-m-1\\
  &= M+m+2-m-1\\
  &=M+1\\
  &> M.
   \nonumber
\end{aligned}
\end{equation}
\noindent Therefore, $k_r=+\infty$.$\hfill\square$\\

In this paper, we always assume that one of them are finite. If there is an $s\in\mathbb Q^{+}$ with $l_s<+\infty$, then for any $r\in\mathbb Q^{+}$, $l_r<+\infty$ by Lemma 4.9 and $k_r=+\infty$ by Lemma 4.10. Similar result for $k_r$ can be obtained.

Now we can define a graded map from $\mathbb Q\longrightarrow \mathbb Z$.\\

\noindent {\bf Lemma 4.11.} Let $A=\bigoplus_{r\in \mathbb Q}{A_rX^r}$ be a graded extension of $V$ in $K[\mathbb Q,\sigma ]$ of Type (e). Assume that there exist $t\in\mathbb Q$ with $O_l(A_t)=W$ and either $l_m$ or $k_m$ is finite for some $m\in\mathbb Q^{+}$. Then there exist a non-zero graded map $f:\mathbb Q\longrightarrow\mathbb Z$ and $\alpha_r$(for any $r\in\mathbb Q)$ with $\alpha_0=1$, $WA_{r}=Wb^{f(r)}\alpha_r$, $W\alpha_r=\alpha_rW^{\sigma(r)}$ for any $r\in\mathbb Q$ and $(W\alpha_r)(W\alpha_s)^{\sigma(r)}=W\alpha_r\alpha_s^{\sigma(r)}=W\alpha_{r+s}$ for any $r,s\in\mathbb Q$.

Furthermore, if $l_m<+\infty$ for some $m\in\mathbb Q^{+}$, then $f=f_{0}^{(-1)}$ or $f=f_d$, $f_{d}^{(1)}$ or $f_{d}^{(-1)}$ for some $d>0$. If $k_m<+\infty$ for some $m\in\mathbb Q^{+}$, then $f=f_{0}^{(1)}$ or $f=f_d$, $f_{d}^{(1)}$ or $f_{d}^{(-1)}$ for some $d<0$.\\

\noindent {\bf Proof.} If $l_m<+\infty$, then by Lemma 4.9, $l_k<+\infty$ for any $k\in\mathbb Q^{+}$. Define $f:\mathbb Q\longrightarrow\mathbb Z$ with $f(0)=0$, for any $k\in\mathbb Q^{+}$, $f(k)=l_k$. Since $l_k\in\mathbb Z$, there exist a $j\in\mathbb N$ such that $g_{\frac{1}{jq}}(k)=l_k$. Set $f(-k)=g_{\frac{1}{jq}}(-k)$. Set $\alpha_0=1$, $\alpha_k=b^{-f(k)}\beta_{k}$ for any $k\in\mathbb Q\backslash\{0\}$. Then $f:\mathbb Q\longrightarrow \mathbb Z$ is a map and for any $k\in\mathbb Q$, $WA_{k}=Wb^{f(k)}\alpha_k$. We claim $f$ is a non-zero graded map. For any $r,s\in\mathbb Q$, $r\neq -s$, $r,s\neq 0$. Then there exist $q_1,j_1,q_2,j_2,q_3,j_3\in\mathbb N$ such that $g_{\frac{1}{j_1q_1}}(\lvert r \rvert)=l_{\lvert r \rvert}$, $g_{\frac{1}{j_2q_2}}(\lvert s \rvert)=l_{\lvert s \rvert}$, $g_{\frac{1}{j_3q_3}}(\lvert r+s \rvert)=l_{\lvert r+s \rvert}$. Let $m=j_1q_1j_2q_2j_3q_3$. Then $\bigoplus_{i\in\mathbb Z}B_{\frac{i}{m}}X^{\frac{i}{m}}$ is a graded extension of $W$ in $K[X^{\frac{1}{m}},X^{-\frac{1}{m}},\sigma(\frac{1}{m})]$ of Type (e). $f(r)=g_{\frac{1}{m}}(r)$, $f(s)=g_{\frac{1}{m}}(s)$, $f(r)+f(s)=g_{\frac{1}{m}}(r)+g_{\frac{1}{m}}(s)\leq g_{\frac{1}{m}}(r+s)=f(r+s)$, $f(r)+f(-r)=g_{\frac{1}{m}}(r)+g_{\frac{1}{m}}(-r)\geq-1$. If $r=-s$, $r,s\neq 0$,  similarly, we can get $f(r)+f(-r)\geq-1$.  Either $r=0$ or $s=0$, clearly $f(r)+f(s)=f(r+s)$. Then $f$ is a graded map. By  our definition, for any $r>0$, $f(r)\geq 0$, $f(-r)<0$. Hence $f=f_{0}^{(-1)}$ or $f=f_d$, $f_{d}^{(1)}$ or $f_{d}^{(-1)}$ for some $d>0$ by Theorem 2.8. Since $\bigoplus_{i\in\mathbb Z}B_{\frac{i}{m}}X^{\frac{i}{m}}$ is a graded extension of Type (e), $WA_r=Wb^{f(r)}\alpha_r$, $W\alpha_r=\alpha_rW^{\sigma(r)}$, $(W\alpha_r)(W\alpha_s)^{\sigma(r)}=W\alpha_r\alpha_s^{\sigma(r)}=W\alpha_{r+s}$ for any $r,s\in\mathbb Q$. If $k_m<+\infty$. Then we consider the graded extension $E=\bigoplus_{r\in\mathbb Q}E_rY^{r}$ of $K[\mathbb Q,\tau]$, we can get similar results.$\hfill\square$\\

\noindent {\bf Lemma 4.12.} Under the assumptions of Lemma 4.11, let $r\in\mathbb Q$, $r\neq 0$, Then $Wb^{f(r)-1}\alpha_r\subset A_r\subset Wb^{f(r)}\alpha_r$ if and only if $W\neq V$ and $f(r)+f(-r)=0$.\\

\noindent {\bf Proof.} Let $r=\frac{n_1}{m_1}$, $m_1\in\mathbb N$, $n_1\in\mathbb Z$, let $lm_1=m$. Then $H_{\frac{1}{m}}=\bigoplus_{i\in\mathbb Z}A_{\frac{i}{m}}X^{\frac{i}{m}}$ is a graded extension of Type (e) with $r=(ln_1)\cdot{\frac{1}{m}}$. So we can get the result by Theorem 2.9.$\hfill\square$\\

By Theorem 2.8, if $f:\mathbb Q\longrightarrow\mathbb Z$ is a non-zero graded map and $f(r)+f(-r)=0$ for some $r\in\mathbb Q^{+}$. Then there exist a smallest positive rational number $l$ with $f(l)+f(-l)=0$. Similar to the proof of [10], we have the following lemma.\\

\noindent {\bf Lemma 4.13.} Let $f:\mathbb Q\longrightarrow\mathbb Z$ be a non-zero graded map. Assume that $f(r)+f(-r)=0$ for some $r\in\mathbb Q^{+}$. Then

 (1) For any $s\in\mathbb Q$, $f(s+r)=f(s)+f(r)$, $f(s-r)=f(s)+f(-r)$.

 (2) Assume that $k$ is the smallest rational number with $f(k)+f(-k)=0$. Then $f(s)+f(-s)=0$ if and only if $s\in k\mathbb Z$.\\

Next theorem is one of the main results of this section.\\

\noindent {\bf Theorem 4.14.} Let $A=\bigoplus_{r\in\mathbb Q}{A_rX^r}$ be an additive subgroup of $K[\mathbb Q,\sigma]$.
$A_0=V$, $W_r=O_(A_r)$, $W=\bigcup_{r\in\mathbb Q}W_r$, $J(W)=Wb^{-1}$ for some $b\in K$. Assume that there is a $t\in\mathbb Q$ with $O_l(A_t)=W$. Also, assume that either $l_m<+\infty$ or $k_m<+\infty$ for some $m\in\mathbb Q^{+}$. Then $A$ is a graded extension of $V$ in $K[\mathbb Q,\sigma]$ of Type (e) if and only if the following properties hold:

 (1) There is a non-zero graded map $f:\mathbb Q\longrightarrow \mathbb Z$. For any $r\in\mathbb Q$, there is an $\alpha_r\in\mathbb Q$ with $\alpha_0=1$. For any $r,s\in\mathbb Q$, $WA_{r}=Wb^{f(r)}\alpha_r$, $W\alpha_r=\alpha_rW^{\sigma(r)}$, $(W\alpha_r)(W\alpha_s)^{\sigma(r)}=W\alpha_r\alpha_s^{\sigma(r)}=W\alpha_{r+s}$.

 (2)(a) If either $W=V$ or $f(r)+f(-r)=-1$ for any $r\neq 0$. Then for any $r\neq 0$, $A_r=Wb^{f(r)}\alpha_r$.

\hspace{1.2em}(b) If $W\neq V$ and there is a $k\in\mathbb Q^{+}$ with $f(k)+f(-k)=0$(assume that $k$ is the smallest positive rational number for this property), then $A_r=Wb^{f(r)}\alpha_r$ for all $r\notin k\mathbb Z$ and $B=\bigoplus_{j\in\mathbb Z}A_{jk}X^{jk}$ is a graded extension of $V$ in $K[X^{k},X^{-k};\sigma(k)]$ with $Wb^{f(jk)-1}\alpha_{jk}\subset A_{jk}\subset Wb^{f(jk)}\alpha_{jk}$ for all $j\in\mathbb Z$.\\

\noindent {\bf Proof.} Suppose that (1) and either (2)(a) or (2)(b) hold. If (2)(a) holds, by (1),

\begin{equation}
\begin{aligned}
   A_{-r}&=Wb^{f(-r)}\alpha_{-r}\\
  &= W\alpha_{-r}(Wb^{f(-r)})^{\sigma(-r)}\\
  &=W(\alpha_{r}^{-1}Wb^{f(-r)})^{\sigma(-r)} \\
  &\supseteq (\alpha_{r}^{-1}b^{-f(r)-1}W)^{\sigma(-r)}\\
  &= (\alpha_{r}^{-1}b^{-f(r)}J(W))^{\sigma(-r)},
   \nonumber
\end{aligned}
\end{equation}

\noindent so $A_r\cup(A_{-r}^{-})^{\sigma(r)}=K$. If (2)(b) holds, for any $r\notin k\mathbb Z$, similarly, $A_r\cup(A_{-r}^{-})^{\sigma(r)}=K$. If $r\in k\mathbb Z$, clearly, $A_r\cup(A_{-r}^{-})^{\sigma(r)}=K$.

Next we will prove that $A$ is a ring. Note that

\begin{equation}
\begin{aligned}
   A_rA_s^{\sigma(r)}&\subseteq (Wb^{f(r)}\alpha_r)(Wb^{f(s)}\alpha_{s})^{\sigma(r)}\\
  &= Wb^{f(r)+f(s)}\alpha_{r}\alpha_{s}^{\sigma(r)}\\
  &=Wb^{f(r)+f(s)}\alpha_{r+s},
   \nonumber
\end{aligned}
\end{equation}

\noindent so if $r+s\notin k\mathbb Z$, then  $A_rA_s^{\sigma(r)}\subseteq A_{r+s}$ follows. Assume that $r+s\in k\mathbb Z$, there are two cases, i.e., either $r,s\in k\mathbb Z$ or $r,s\notin k\mathbb Z$. If $r,s\in k\mathbb Z$, then $A_rA_s^{\sigma(r)}\subseteq A_{r+s}$, since $B$ is a graded extension of $V$ in $K[X^{k},X^{-k};\sigma(k)]$. If $r,s\notin k\mathbb Z$, then $s=mk-r$ for some $m\in\mathbb Z$. So we have

\begin{equation}
\begin{aligned}
   f(r)+f(s)&= f(r)+f(-r+mk)\\
  &= f(r)+f(-r)+f(mk)\\
  &=-1+f(mk)\\
  &=-1+f(r+s)
   \nonumber
\end{aligned}
\end{equation}

\noindent by Lemma 4.13. Thus $A_rA_s^{\sigma(r)}\subseteq Wb^{f(r)+f(s)}\alpha_{r+s}=Wb^{f(r+s)-1}\alpha_{r+s}\subset A_{r+s}$. Hence $A$ is a graded extension of $V$ in $K[\mathbb Q,\sigma]$ by Lemma 2.5.

 Conversely, suppose that $A=\bigoplus_{r\in\mathbb Q}A_rX^{r}$ is a graded extension of $V$ in $K[\mathbb Q,\sigma]$. Then (1) holds by Lemma 4.11, (2)(a) holds by Lemma 4.12 and (2)(b) holds by Lemmas 4.12 and 4.13.$\hfill\square$\\

 At last, we will study extensions of Type (h). We always assume that $W$ is a total valuation ring of $K$ with $W\supset V$, $J(W)^2=J(W)$.\\

 \noindent {\bf Lemma 4.15.} Let $A_r$ be a left $V$-ideal of $K$ and $O_l(A_r)\subset W$. Then $WA_{r}$ is a princial left $W$-ideal.\\

 \noindent {\bf Proof.} Since $O_l(A_r)\subset W$, $J(W)A_r\subset A_r\subset WA_r$. If $WA_r$ is not a principal left $W$-ideal, then $WA_r=J(W)WA_r=J(W)A_r$, a contradiction.$\hfill\square$\\

 \noindent {\bf Lemma 4.16.} Let $A_r$ be a left $V$-ideal of $K$, $O_l(A_r)=W_r\subseteq W$, $J(W)^2=J(W)$, $^\ast(J(W)A_r)=J(W)A_r$. Then $A_r=WA_r$ is not a principal left $W$-ideal and $^\ast A_r=A_r$.\\

 \noindent {\bf Proof.} Suppose that $WA_r=Wc$ for some $c\in K$. Then
 \begin{center}
 $^\ast(J(W)A_r)=^\ast(J(W)c)=Wc\neq J(W)A_r$,
 \end{center}
 \noindent a contradiction. So $WA_r$ is not a principal left $W$-ideal. Then $O_l(A_r)=W$ by Lemma 4.15. Hence $A_r=WA_r$ is not principal. Therefore $A_r=J(W)A_r$ and $^\ast A_r=A_r$.$\hfill\square$\\

 \noindent {\bf Lemma 4.17.} Let $A=\bigoplus_{r\in \mathbb Q}{A_rX^r}$ be a graded extension of $V$ in $K[\mathbb Q,\sigma ]$ of Type (h), $W_r=O_l(A_r)$, $W=\bigcup_{r\in\mathbb Q}W_r$. Assume that there is a $t\in\mathbb Q$ with  $O_l(A_t)=W$. Then there exist an $l\in\mathbb N$ such that $H_{\frac{1}{l}}=\bigoplus_{i\in \mathbb Z}{A_{\frac{i}{l}}X^{\frac{i}{l}}}$ is a graded extension of $V$ in $K[X^{\frac{1}{l}},X^{-\frac{1}{l}};\sigma(\frac{1}{l})]$ of Type (h) with $O_l(A_{\frac{1}{l}})=W$.\\

 \noindent {\bf Proof.} We can assume that $t\neq 0$. Let $t=\frac{n_1}{m_1}$, $m_1\in\mathbb N$, $n_1\in\mathbb Z$. Since $A$ is graded extension of Type (h), there is an $s\in\mathbb Q^{+}$ with $H_s=\bigoplus_{i\in\mathbb Z}A_{is}X^{is}$ is a graded extension of $V$ in $K[X^s,X^{-s};\sigma(s)]$ of Type (II). Let $s=\frac{n_2}{m_2}$, $m_2,n_2\in\mathbb N$. Let $l=m_1m_2$. Then $H_{\frac{1}{l}}$ is a graded extension of $V$ in $K[X^{\frac{1}{l}},X^{\frac{1}{l}};\sigma({\frac{1}{l}})]$ of Type (h) with $O_l(A_{\frac{1}{l}})=W$.$\hfill\square$\\

\noindent {\bf Lemma 4.18.}  Under the assumptions of Lemma 4.17. Let $r\in\mathbb Q$, $r\neq 0$. Assume $^\ast(J(W)A_r)=Wc$ for some $c\in K$, then $^\ast(J(W)A_r)=Wc=cW^{\sigma(r)}$.\\

\noindent {\bf Proof.} Since $A$ is  graded extension of Type (h), there exist an $l\in\mathbb N$ such that $H_{\frac{1}{l}}=\bigoplus_{i\in \mathbb Z}{A_{\frac{i}{l}}X^{\frac{i}{l}}}$ is a graded extension of $V$ in $K[X^{\frac{1}{l}},X^{-\frac{1}{l}};\sigma(\frac{1}{l})]$ of Type (h) with $O_l(A_{\frac{1}{l}})=W$. Let $r=\frac{n_1}{m_1}$, $m_1\in\mathbb N$, $n_1\in\mathbb Z$, $p=lm_1$. Then $H_{\frac{1}{p}}$ is a graded extension of $V$ in $K[X^{\frac{1}{p}},X^{-\frac{1}{p}};\sigma(\frac{1}{p})]$ of Type (h) with $O_l(A_{\frac{1}{p}})=W$. $r=\frac{n_1}{m_1}=(n_1l)\cdot{\frac{1}{p}}$. Then $^\ast(J(W)A_r)=Wc=cW^{\sigma(r)}$ can be get by [11, Lemma 1.18].$\hfill\square$\\

 \noindent {\bf Lemma 4.19.} Under the assumptions of Lemma 4.17. Let
\begin{center}
$H=\{r\in\mathbb Q\mid ^\ast(J(W)A_r)=Wc_r$
for some $c_r\in K\}$.
\end{center}
\noindent Then $H$ is an additive subgroup of $\mathbb Q$.\\

\noindent {\bf Proof.} Obviously, $0\in H$. For any $r,s\in H$, let $r=\frac{n_1}{m_1}$, $s=\frac{n_2}{m_2}$, $m_1,m_2\in\mathbb N$, $n_1,n_2\in\mathbb Z$. Set $p=m_1m_2l$. Then $H_{\frac{1}{p}}$ is a graded extension of $V$ in $K[X^{\frac{1}{p}},X^{-\frac{1}{p}};\sigma(\frac{1}{p})]$ of Type (h). $r=n_1/m_1=(n_1m_2l)\cdot{\frac{1}{p}}$, $s=n_2/m_2=(n_2m_1l)\cdot{\frac{1}{p}}$. Hence $-r\in H$, $r+s\in H$ by [11, Lemma 1.16].$\hfill\square$\\

\noindent {\bf Corollary 4.20.} Under the assumptions and notations of Lemma 4.19, for any $r\notin H$, $^\ast A_r=A_r=WA_r$ is not principal left $W$-ideal.\\

\noindent {\bf Proof.} The result can be obtained by Lemmas 4.16 and 4.19.$\hfill\square$\\

\noindent {\bf Lemma 4.21.} Under the assumptions and notations of Lemma 4.19, we have $H\neq \{0\}$, $H\neq\mathbb Q$.\\

\noindent {\bf Proof.} Since $H_{\frac{1}{l}}=\bigoplus_{i\in \mathbb Z}{A_{\frac{i}{l}}X^{\frac{i}{l}}}$ is a graded extension of $V$ in $K[X^{\frac{1}{l}},X^{-\frac{1}{l}};\sigma(\frac{1}{l})]$ of Type (h) with $O_l(A_{\frac{1}{l}})=W$. So there is an $m\in\mathbb N$ with $\frac{m}{l}\in H$. Hence $H\neq \{0\}$. $\frac{1}{l}\notin H$ implies $H\neq\mathbb Q$.$\hfill\square$\\

Next theorem is another main result of this section.\\

\noindent {\bf Theorem 4.22.} Let $A=\bigoplus_{r\in \mathbb Q}{A_rX^r}$ be an additive subgroup of $K[\mathbb Q,\sigma]$, $A_0=V$, $W_r=O_l(A_r)$, $W=\bigcup_{r\in\mathbb Q}W_r$, $J(W)^2=J(W)$. Assume that there is a $t\in\mathbb Q$ with $O_l(A_t)=W$. Let $H=\{r\in\mathbb Q\mid ^\ast(J(W)A_r)=Wc_r$ for some $c_r\in K\}$. Then $A$ is a graded extension of $V$ in $K[\mathbb Q,\sigma]$ of Type (h) if and only if for any $r,s\in H$, $^\ast(J(W)A_s)=Wc_s$ for some $c_s\in K$, $Wc_s=c_sW^{\sigma(s)}$, $Wc_{-s}=W(c_s^{-1})^{\sigma(-s)}$, $(Wc_r)(Wc_s)^{\sigma(r)}=Wc_{r+s}$ and the following properties hold:

(1) $H$ is an additive non-zero proper subgroup of $\mathbb Q$ with $B=\bigoplus_{s\in H}A_sX^{s}$ is a graded extension of $V$ in $K[H,\sigma\mid_{H}]$.

 (2) For any $s\notin H$, we have $^\ast A_s=A_s=WA_s$ is not a principal left $W$-ideal and $O_r(A_s)=W^{\sigma(s)}$.

 (3) $M=\bigoplus_{r\in\mathbb Q}M_rX^{r}=V\bigoplus(\bigoplus_{-s\in H^{-}}
Wc_{-s}X^{-s})\bigoplus(\bigoplus_{s\in H^{+}}J(W)c_sX^{s})\\
\bigoplus(\bigoplus_{r\notin H}A_rX^{r})$ is a graded extension of $V$ in $K[\mathbb Q,\sigma]$, where $H^{+}=\{r\in H\mid r>0\}$, $H^{-}=\{r\in H\mid r<0\}$.\\

\noindent {\bf Proof.} Let $A$ be a graded extension of $V$ in $K[\mathbb Q,\sigma]$. (1), (2) hold by corollary 4.20 and Lemmas 4.18, 4.19, 4.21. Let $r_1,r_2\in\mathbb Q$. If $r_1+r_2\in H$, then $r_1,r_2\in H$ or $r_1,r_2\notin H$. If $r_1,r_2\notin H$, then $M_{r_1}(M_{r_2})^{\sigma(r_1)}=A_{r_1}(A_{r_2})^{\sigma(r_1)}\subseteq J(W)A_{r_1+r_2}\subseteq M_{r_1+r_2}$. Let $r_1\in H$, $r_2\in H$. Either $r_1=0$ or $r_2=0$, then $M_{r_1}(M_{r_2})^{\sigma(r_1)}=M_{r_1+r_2}$. Assume $r_1\neq 0$, $r_2\neq 0$. If either $r_1>0$ or $r_2>0$, then $M_{r_1}(M_{r_2})^{\sigma(r_1)}=(J(W)A_{r_1})(J(W)A_{r_2})^{\sigma(r_1)}\subseteq J(W)A_{r_1+r_2}\subseteq M_{r_1+r_2}$. If $r_1,r_2<0$, then $M_{r_1}(M_{r_2})^{\sigma(r_1)}=(Wc_{r_1})(Wc_{r_2})^{\sigma(r_1)}=Wc_{r_1+r_2}= M_{r_1+r_2}$. If $r_1+r_2\notin H$, then $r_1\notin H$ or $r_2\notin H$. Hence $M_{r_1}(M_{r_2})^{\sigma(r_1)}=J(W)A_{r_1}(A_{r_2})^{\sigma(r_1)}\subseteq J(W)A_{r_1+r_2}=A_{r_1+r_2}=M_{r_1+r_2}$. Obviously,
\begin{center}
$M_r\bigcup (M_{-r}^{-})^{\sigma(r)}=K$ for any $r\in\mathbb Q$.
\end{center}
\noindent Hence $M$ is a graded extension of $V$ in $K[\mathbb Q,\sigma]$.

Conversely, suppose that (1), (2), (3) hold. Then we have $A_r\bigcup(A_{-r}^{-})^{\sigma(r)}=K$ for any $r\in\mathbb Q$. For any $r_1,r_2\in\mathbb Q$, let $r_1=\frac{n_1}{m_1}$, $r_2=\frac{n_2}{m_2}$, $m_1,m_2\in\mathbb N$, $n_1,n_2\in\mathbb Z$. Since $H\neq \mathbb Q$, there exist $r_3\notin H$, $^\ast A_{r_3}=A_{r_3}=WA_{r_3}$ is not a principal left $W$-ideal, $H\neq \{0\}$ implies there exist an $s\in H^{+}$. Let $r_3=\frac{n_3}{m_3}$, $s=\frac{n_4}{m_4}$, $m_3,m_4\in\mathbb N$, $n_3,n_4\in\mathbb Z$. Let $t=\frac{n_5}{m_5}$, $m_5\in\mathbb N$, $n_5\in\mathbb Z$. Let $p=m_1m_2m_3m_4m_5$. Then $r_1=(n_1m_2m_3m_4m_5)\cdot{\frac{1}{p}}$, $r_2=(n_2m_1m_3m_4m_5)\cdot{\frac{1}{p}}$. Since $M$ is a graded extension of $V$ in $K[\mathbb Q,\sigma]$, $G_r=\bigoplus_{i\in\mathbb Z}M_{\frac{i}{p}}X^{\frac{i}{p}}$ is a graded extension of Type (h). By Theorem 2.10, $H_r=\bigoplus_{i\in\mathbb Z}A_{\frac{i}{p}}X^{\frac{i}{p}}$ is a graded extension of Type (h). Hence $A_{r_1}A_{r_2}^{\sigma(r_1)}\subseteq A_{r_1+r_2}$. Therefore, $A$ is a graded extension of $V$ in $K[\mathbb Q,\sigma]$.$\hfill\square$\\

\noindent {\bf Remark 4.23.} Let $A=\bigoplus_{r\in\mathbb Q}{A_r}X^{r}$ be a graded extension of $V$ in $K[\mathbb Q,\sigma]$ of Type (II), $W_r=O_l(A_r)$ for any $r\in\mathbb Q$, $W=\bigcup_{r\in\mathbb Q}W_r$. We have two conjectures:

 (1) There is always a $t\in\mathbb Q$ with $O_l(A_t)=W$.

 (2) Let $A$ be a graded extension of Type (e). Then either $l_m<+\infty$ or $k_m<+\infty$ for some $m\in\mathbb Q^{+}$, where $l_m$ and $k_m$ are defined in this section.

If both of these conjectures are true, then we can describe the structure of Type (II) completely.
Unfortunately, we can not find methods to prove them.\\[1.0cm]
{\bf 5. Examples}\\

\noindent In this section, we will provide concrete examples of graded extensions of $V$ in $K[\mathbb Q,\sigma]$ for illustrating the classification.

Let $W$ be an overring of $V$ with $J(W)^2=J(W)$, $J(V)=b^{-1}V$ for some $b\in K$ and $\sigma(r)=1_{K}$ for any $r\in\mathbb Q$, Then the following are trivial.\\

\noindent {\bf Example 5.1.} (1) $A=\bigoplus_{r\in\mathbb Q}VX^{r}$ is a graded extension of $V$ in $K\mathbb Q$ of Type (a).

 (2) $A=(\bigoplus_{r\in\mathbb Q^{+}}J(W)X^{r})\bigoplus V\bigoplus (\bigoplus_{r\in\mathbb Q^{-}}WX^{r})$ is a graded extension of $V$ in $K\mathbb Q$ of Type (f).

 (3) $A=\bigoplus_{r\in\mathbb Q}Vb^{[r]}X^{r}$ is a graded extension of $V$ in $K\mathbb Q$ of Type (e).\\

Let $F_0[Y^r\mid r\in\mathbb Q]$ be a commutative domain over a field $F_0$ in indeterminate $Y$ with $Y^s\cdot Y^t=Y^{s+t}$ and let $K=F_0(Y^r\mid r\in\mathbb Q)$ be its quotient field. We define a map $v$ from $K\backslash\{0\}$ to $\mathbb Q$ as follows;

\begin{center}
$v(\sum a_{r_i}Y^{r_i})=min\{r_i\mid a_{r_i}\neq 0\}$
\end{center}

\noindent and $v(\alpha\beta^{-1})=v(\alpha)-v(\beta)$ for any $\alpha,\beta\in F_0[Y^r\mid r\in\mathbb Q]$ with $\beta\neq 0$. Then $v$ is a valuation of $K$. Let $V=\{a\mid v(a)\geq 0, a\in K\backslash \{0\}\}\bigcup \{0\}$. Then $V$ is a valuation ring of $K$. Let $\sigma(r)=1_{K}$ for any $r\in\mathbb Q$. Let $\pi$ be a positive real number but not a rational number. Then we have the following example.\\

\noindent {\bf Example 5.2.} For any $r>0$, let $A_r=\bigcup\{VY^{-s}\mid s<r\pi,s\in\mathbb Q\}$ and $A_{-r}=\bigcup\{VY^{s}\mid s>r\pi,s\in\mathbb Q\}$. Then $O_l(A_r)=V$, $J(V)^2=J(V)$ and $A_s^\ast=A_s$ is not a principal $V$-ideal for any $s\neq 0$. Let $A_0=V$. Set $A=\bigoplus_{r\in\mathbb Q}A_rX^r$. Then $A$ is a graded extension of $V$ in $K\mathbb Q$ of Type (g).\\

The following examples are similar to [11, Example 2.2].

Let $H$ be a non-zero proper subgroup of $\mathbb Q$. Let $F_0[Y_r^s\mid r\in\mathbb Q,s\in \mathbb Q+H\pi]$ be a commutative domain over a field $F_0$ in indeterminates $Y_r$ with $Y_r^{s_1} \cdot Y_r^{s_2}=Y_r^{s_1+s_2}$ and let $F=F_0(Y_r^s\mid r\in\mathbb Q,s\in\mathbb Q+H\pi)$ be its quotient field. We define a group homomorphism $\sigma:\mathbb Q\longrightarrow Aut(F)$ which satisfies the following; for any $t\in\mathbb Q$, $\sigma(t)(a)=a$ for all $a\in F_0$ and $(Y_r^{s})^{\sigma(t)}=Y_{r-t}^s$ for any $r\in\mathbb Q$, $s\in\mathbb Q+H\pi$.

Furthermore, let $G=\bigoplus_{s\in\mathbb Q}{\mathbb R_s}$, the direct sum of $\mathbb R_s$ with $\mathbb R_s=\mathbb R$. We define a total order on $G$ as follows; for any $\gamma_1=(a_r)_{r\in\mathbb Q}$, $\gamma_2=(b_s)_{s\in\mathbb Q}\in G$, $\gamma_2<\gamma_1$ if and only if there exist a $k\in\mathbb Q$ with $a_r=b_r$ for all $r<k$ and $b_k<a_k$.
We define a map $v$ from $F\backslash\{0\}$ to $G$ as follows; $v(a)=0$ for any non-zero $a\in F_0$ and for any non-zero homogeneous element $\alpha=Y_{r_1}^{s_1}\cdot\cdot\cdot Y_{r_n}^{s_n}(r_1<r_2<\cdot\cdot\cdot<r_n)$, $v(\alpha)=(s_k)_{k\in\mathbb Q}$ with the $r_j$-component of $v(\alpha)$ is $s_j(1\leq j\leq n)$ and the other components of it are all zeroes.

 Let $\beta=\beta_1+\cdot\cdot\cdot+\beta_m$ be any element in $F_0[Y_r^s\mid r\in\mathbb Q,s\in\mathbb Q+H\pi]$, where $\beta_i$ are non-zero homogenous elements. Then we define $v(\beta)=min\{v(\beta_{i})\mid 1\leq i\leq m\}$. As usual, we can extend the map $v$ to $F\backslash\{0\}$, which is a valuation of $F$. We denote by $V_0$ the valuation ring of $F$ determined by $v$, that is, $V_0=\{\alpha\beta^{-1}\mid v(\alpha\beta^{-1})=v(\alpha)-v(\beta)\geq 0$, where $\alpha,\beta\in F_0[Y_r^s\mid r\in\mathbb Q,s\in\mathbb Q+H\pi]$ with $\alpha,\beta\neq 0\}\bigcup \{0\}$.

Note that $\sigma(t)(V_0)=V_0$ for any $t\in\mathbb Q$, since $\sigma(t)$ is just shifting and that, for any $\alpha\beta^{-1}\in F$, $V_0\alpha\beta^{-1}=V_0Y_{r_1}^{s_1}\cdot\cdot\cdot Y_{r_n}^{s_n}$ for some homogenous element $Y_{r_1}^{s_1}\cdot\cdot\cdot Y_{r_n}^{s_n}$ by the construction of $v$. We set $U_r=\bigcup\{V_0Y_r^s\mid s\in \mathbb Q+H\pi\}$, an overring of $V_0$. Then $(U_r)^{\sigma(t)}=U_{r-t}\supset U_r$ for any $t>0$.

Let $S=F[\mathbb Q,\sigma]$ be the skew group ring of $\mathbb Q$ over $F$. By [8, Theorem 2.2], $S$ has a quotient ring $K$. Let $S=\{\sum a_{s_i}Z^{s_i}\mid a_{s_i}\in F,{s_i}\in\mathbb Q\}$. For any $r\in\mathbb Q$, $f=\sum a_{s_i}Z^{s_i}\in S$, let $\sigma(r)(f)=\sum (a_{s_i})^{\sigma(r)}Z^{s_i}$, for any $\alpha=fg^{-1}\in K$, $\sigma(r)(fg^{-1})=f^{\sigma(r)}(g^{\sigma(r)})^{-1}$. Then $\sigma(r)$ is an automorphism of $K$. Then $\sigma:\mathbb Q\longrightarrow Aut(K)$ is a group homomorphism.

Let $R=\{\sum a_{s_i}Z^{s_i}\mid a_{s_i}\in F,s_i\in\mathbb Q,s_i\geq 0\}$. Then $R$ is a domain. Let $I$ be the ideal of $R$ generated by $Z^s(s>0)$. Then $I$ is a completely prime ideal of $R$ and $I$ is localizable. Let $R_I$ be the localization of $R$ at the ideal $I$. Then $R_{I}$ is a total valuation ring.

 Similar to [11, Example 2.2], we define the map

\begin{center}
$\varphi:R_I\longrightarrow F$
\end{center}

\noindent by $\varphi(fg^{-1})=f(0)g(0)^{-1}$, where $f,g\in R$, $g(0)\neq 0$. We let

\begin{center}
$V=\varphi^{-1}(V_0)$ and $W_r=\varphi^{-1}(U_r)$,
\end{center}

\noindent the complete inverse image of $V_0$ and $U_r$ by $\varphi$ respectively for any $r\in\mathbb Q$. Similar to [6, Proposition 1.6], we can get $V$ and $W_r$ are all total valuation rings of $K$. Furthermore, we have the following properties which are easily proved by the definitions:

(1) $V^{\sigma(s)}=V$  for any $s\in\mathbb Q$.

(2) $(W_r)^{\sigma(s)}=W_{r-s}\supset W_{r}$ for any $r\in\mathbb Q$, $s\in\mathbb Q^{+}$.

(3) $Y_r^sV=VY_r^s$ and $Y_r^sW_t=W_tY_r^s$ for any $r,t\in\mathbb Q$, $s\in \mathbb Q+H\pi$.

(4) $J(V)^2=J(V)$ and $J(W_r)^2=J(W_r)$ for any $r\in\mathbb Q$.

(5) $Z^{r}W_s=W_{s-r}Z^{r}$ for any $r,s\in\mathbb Q$.\\

\noindent {\bf Example 5.3.} (1) Let $A=V\bigoplus(\bigoplus_{r\in\mathbb Q^{+}}VZ^rX^{r})\bigoplus(\bigoplus_{r\in\mathbb Q^{+}}J(V)Z^{-r}X^{-r})$. Then $A$ is a graded extension of $V$ in $K[\mathbb Q,\sigma]$ of Type (d).

 (2)  Let $A=V\bigoplus(\bigoplus_{r\in\mathbb Q^{+}}J(V)Z^rX^{r})\bigoplus(\bigoplus_{r\in\mathbb Q^{+}}VZ^{-r}X^{-r})$. Then $A$ is a graded extension of $V$ in $K[\mathbb Q,\sigma]$ of Type (f).
\\

\noindent {\bf Example 5.4.} Let

\begin{center}
$A=\bigoplus_{r\in\mathbb Q}A_{r}X^{r}=V\bigoplus(\bigoplus_{r\in\mathbb Q^{+}}W_0Z^{-2r}X^{r})\bigoplus(\bigoplus_{r\in\mathbb Q^{+}}J(W_{-r})Z^{2r}X^{-r})$.
\end{center}
\noindent For any $r>0$, $W_0Z^{-2r}\supset W_rZ^{-2r}=Z^{-2r}(W_0)^{\sigma(r)}$. We can check easily that $A$ is a graded extension of $K[\mathbb Q,\sigma]$ of Type (b).\\

\noindent {\bf Example 5.5.} Let

\begin{center}
$A=\bigoplus_{r\in\mathbb Q}A_rX^r=V\bigoplus(\bigoplus_{r\in\mathbb Q^{+}}W_{-2r}Z^{r}X^{r})\bigoplus(\bigoplus_{r\in\mathbb Q^{+}}J(W_0)Z^{-r}X^{-r})$.
\end{center}

\noindent For any $r>0$, $W_{-2r}Z^{r}\subset W_{-4r}Z^{r}=Z^{r}(W_{-2r})^{\sigma(r)}$. We can check that $A$ is a graded extension of $K[\mathbb Q,\sigma]$ of Type (c).\\

\noindent {\bf Remark.} In Examples 5.4 and 5.5, set $W=\bigcup_{r\in\mathbb Q}O_l(A_r)$. Then for any $r\in\mathbb Q$,  $O_l(A_r)\subset W$.\\

\noindent {\bf Example 5.6.} Let $r\in\mathbb Q$, $r>0$. If $r\notin H$, set
\begin{center}
$A_r=\bigcup\{VY_0^{-s}Z^{-r}\mid s<r\pi, s\in \mathbb Q+H\pi\}$,
$A_{-r}=\bigcup\{Z^{r}Y_r^{s}V\mid s>r\pi, s\in \mathbb Q+H\pi\}$.
\end{center}
 If $r\in H$, set
 \begin{center}
 $A_r=VY_0^{-r\pi}Z^{-r}$, $A_{-r}=Z^{r}Y_r^{r\pi}J(V)$, $A_0=V$.
 \end{center}
  Then similar to [11, Example 2.2], we have

 (1) $O_l(A_r)=\bigcup_{r>0}W_r$ if $r\notin H$. $O_l(A_r)=V$ if $r\in H$.

 (2) $W=\bigcup_{r\in\mathbb Q}O_l(A_r)=\bigcup_{r>0}W_r$, $J(W)=\bigcup_{s>0}Y_{0}^sW$, $J(W)^2=J(W)$.

(3) For any $r\in H$, $J(W)A_r=J(W)Y_0^{-r\pi}Z^{-r}$.

 (4) If $r\notin H$, $^\ast A_r=A_r$ is not a finite generated $W$-ideal.
 
 Hence we can check that $A=\bigoplus_{r\in\mathbb Q}A_rX^r$ is a graded extension of $V$ in $K[\mathbb Q,\sigma]$ of Type (h).\\[1.0cm]
{\bf Acknowledgments}\\

\noindent G. M. Xie is supported by National Nature Science Foundation of China (11161005) and Guangxi Science Foundation (0991020).\\[1.0cm]

\noindent {\bf References}\\

[1] H. H. Brungs, H. Marubayashi and E. Osmanagic, Gauss extensions and total graded subrings for crossed product algebras, {\textit {J. Algebra}} {\bf 316} (2007) 189-205.

[2] H. H. Brungs and M. Schr\"{o}der, Valuation rings in Ore extensions, {\textit {J. Algebra}} {\bf 235} (2001) 665-680.

[3] H. H. Brungs and G. T\"{o}rner, Extensions of chain rings, {\textit {Math. Z.}} {\bf 185} (1984) 93-104.

[4] H. Marubayashi, H. Miyamoto and A. Ueda,  {\textit {Non-commutative Valuation Rings and Semihereditary Orders,}} (K-Monographs in Math. {\bf 3},  Kluwer Academic Publishers, Boston 1997).

[5] G. Xie, Y. Chen, H. Marubayashi and Y. Wang, A new classification of graded extensions in a skew Laurent polynomial ring. {\textit {Far East J. Math. Sci.}} {\bf 40} (2010) 37-44.

[6] G. Xie, S. Kobayashi, H. Marubayashi, N. Popescu and C. Vraciu, Non-commutative valuation
rings of the quotient Artinian ring of a skew polynomial ring, {\textit {Algebr. Represent. Theory.}} {\bf 8} (2005) 57-68.

[7] G. Xie, S. Kobayashi, H. Marubayashi and H. Komatsu, Non-commuta\\
tive valuation rings of $K(X;\sigma,\delta)$
 over a division ring $K$, {\textit {J. Math. Soc. Japan.}} {\bf 56} (2004)  737-752.

[8] G. Xie, J. Liang and M. Wang, Quotient skew fields of skew group rings of torsion free additive groups over a skew field(preprint).

[9] G. Xie, F. Liu and C. Wei, Graded maps over $\mathbb Q$ and graded extension in $K[\mathbb Q,\sigma]$ of Type (e),  {\textit {Journal of Guangxi Normal University(Natural Science Edition).}} {\bf 28(2)} (2010) 42-46.

[10] G. Xie and H. Marubayashi, A classifification of graded extensions in a skew Laurent polynomial ring, {\textit {J. Math. Soc. Japan.}} {\bf 60} (2008) 423-443.

[11] G. Xie and H. Marubayashi, A classification of graded extensions in a skew Laurent polynomial ring, II, {\textit {J. Math. Soc. Japan.}} {\bf 61(4)} (2009) 1111-1130.

\end{document}